\documentclass[12pt,a4paper]{article}
\usepackage{mathtools}
\usepackage[a4paper,margin=2.2cm]{geometry}
\usepackage[USenglish]{babel}
\DeclareFontFamily{U}{bbold}{}
\DeclareFontShape{U}{bbold}{m}{n}{<-> bbold10}{}

\usepackage{amsthm}
\usepackage{ dsfont }
\let\oldmathbb\mathbb
\let\mathbb\oldmathbb

\usepackage{amssymb}

\usepackage{algorithm}
\usepackage{appendix}
\usepackage[square,authoryear]{natbib}
\usepackage[noend]{algpseudocode}
\usepackage{pgfplots}
\pgfplotsset{compat = newest}
\usetikzlibrary{automata,arrows,positioning,calc}
\usepackage{graphicx}
\usepackage[edges]{forest}
\usepackage{cancel}
\usepackage{textcmds} 

\usepackage{booktabs}
\usepackage{float}

\usepackage{diagbox}

\usepackage{pdfpages} 
\usepackage{lineno}
\usepackage{amsmath}  

\newcommand*\linenomathpatch[1]{%
  \cspreto{#1}{\linenomath}%
  \cspreto{#1*}{\linenomath}%
  \csappto{end#1}{\endlinenomath}%
  \csappto{end#1*}{\endlinenomath}%
}

\linenomathpatch{equation}
\linenomathpatch{gather}
\linenomathpatch{multline}
\linenomathpatch{align}
\linenomathpatch{alignat}
\linenomathpatch{flalign}

\usepackage{subcaption}
\usepackage{caption}
\usepackage{microtype}
\usepackage{multicol}
\usepackage[T1]{fontenc}
\usepackage[]{xcolor}
\usepackage[pdfpagelabels,colorlinks=true,allcolors=black,pdfusetitle]{hyperref}
\usepackage{soul}

\usepackage{float}

\definecolor{ptbblau}{rgb}{0.0,0.608,0.808}
\definecolor{ptbblaubeamer}{rgb}{0.0,0.451,0.667}

\usepackage[right]{eurosym}
\usepackage{tikz}
\usepackage{url}
\usepackage[utf8]{inputenc}

\usepackage{enumitem}

\usepackage[normalem]{ulem}

\newcommand{\R}{\mathbb{R}}

\newcommand{\be}{\begin{eqnarray*}}
\newcommand{\ee}{\end{eqnarray*}}
\newcommand{\ben}{\begin{eqnarray}}
\newcommand{\een}{\end{eqnarray}}
\newcommand{\bi}{\begin{itemize}}
\newcommand{\ei}{\end{itemize}}

\newcommand{\cF}{\mathcal{F}}

\newcommand{\cT}{\mathcal{T}}

\newcommand{\sI}{\text{\usefont{U}{bbold}{m}{n}1}}

\newcommand{\N}{\mathbb{N}}
\newcommand{\IR}{\mathbb{R}}
\newcommand{\IN}{\mathbb{N}}

\newcommand{\IE}{\mathbb{E}}

\newcommand{\td}{\text{d}}

\newcommand{\la}{\lambda}

\newcommand{\si}{{\sigma}}

\newcommand{\weg}[1]{}

\newtheorem{thm}{Theorem}
\newtheorem{prop}[thm]{Proposition}

\newtheorem{lemm}[thm]{Lemma}
\theoremstyle{definition}

\theoremstyle{definition}
\newtheorem{exa}{Example}
\newtheorem{rem}{Remark}

\newcommand{\negpar}[1][-1em]{%
  \ifvmode\else\par\fi
  {\parindent=#1\leavevmode}\ignorespaces
}
\makeatletter
\newcommand{\mylabel}[2]{#2\def\@currentlabel{#2}\label{#1}}
\title{Learning While Searching in the Continuous-Time Parking Problem}

\author{Stefan Ankirchner \thanks{Stefan Ankirchner, Institute for Mathematics, University of Jena, Inselplatz 5, 07743 Jena, Germany. \emph{Email:} s.ankirchner@uni-jena.de.} 
\and Hugalf Bernburg \thanks{Hugalf Bernburg, Institute for Mathematics, University of Jena, Inselplatz 5, 07743 Jena, Germany. \emph{Email:} h.bernburg@uni-jena.de.} 
\and Maximilian Philipp Thiel \thanks{Maximilian Philipp Thiel, Institute for Mathematics, University of Jena, Inselplatz 5, 07743 Jena, Germany. \emph{Email:} maximilian.thiel@uni-jena.de.}  }
\begin{document}
\maketitle
\begin{abstract}
We consider a continuous-time version of the parking problem, a classical example of a stopping problem.
We assume that a stopping agent does not know the distribution of free parking lots, modeled as arrivals of a homogeneous Poisson process.
We allow for rather general cost functions, including the case where the agent represents a group searching for several parking lots. 
We characterize the optimal stopping rule under Bayesian inference in terms of the solution function of the associated HJB equation. We prove existence of a bounded solution of the HJB equation via a continuation method and verify that the solution coincides with the value function. Finally, we demonstrate that in general the optimal stopping rule is not a simple One-Step-Look-Ahead rule. 
\end{abstract}
\begin{center}\footnotesize
  \begin{tabular}{r@{ : }p{10cm}}
      {\itshape 2020 MSC} & Primary: 60G40, 93E35; secondary: 62L05, 68Q32.\\
      {\itshape Keywords} &  optimal stopping, parking problem, reinforcement learning, regret. 
  \end{tabular}
\end{center}

\section*{Introduction}

The parking problem is a classical example of a stopping problem. In the standard continuous-time version an agent is assumed to be driving along a street, where the arrival of free parking lots is modeled with a homogeneous Poisson process. The agent seeks to park as close to some target location as possible. Moreover, it is assumed that the agent cannot make a U-turn and hence, if a free lot is not taken, it cannot be taken anymore later. The decision rules of the agent are modeled as stopping times taking values in the set of jump times of the Poisson process. The parking problem consists in finding the stopping time that minimizes the expected distance of the taken lot to the target location. 

In the classical parking problem the agent is assumed to know the intensity of the Poisson process modeling the arrivals of free parking lots. In this paper we consider the case where the driver does \textit{not} know the jump intensity. How can one characterize optimal decision rules in this case? 
It is important to note that the agent can continuously learn about the jump intensity while driving along the street. In order to determine an optimal stopping rule that continuously takes into account all of the available information we introduce the following Bayesian model. 

We assume that the street starts at $S \in (-\infty, 0)$ and extends to infinity. Let $\Lambda$ be a gamma distributed random variable and $N = (N_t)_{t \in [S, \infty)}$ a homogeneous Poisson process with random jump intensity $\Lambda$, both defined on the same probability space $(\Omega,\cF,P)$.
Let $\cF^N=(\cF^N_t)_{t\in[S,\infty)}$ be the natural filtration of $N$.
Let $\cT$ be the family of $\cF^N$-stopping times $\tau$ with values at jump times of $N$ only.
Note that the filtration $\cF^N$ represents the continuous information flow available for inference of the realization of $\Lambda$. Thus, a stopping time in $\cT$ can make, at any time, use of all the information available. 
Let $c: [S,\infty)\times \IN_0 \to [0,\infty)$ be a cost function.  
If the agent stops at $\tau$, then the expected costs are given by $\IE[c(\tau,N_\tau)]$. We assume that the agent seeks to minimize the expected costs; in other words, she aims at finding $\tau^* \in \cT$ such that 
\begin{align}\label{opti pb}
    \IE[c(\tau^*,N_{\tau^*})]=\min_{\tau\in \cT} \IE[c(\tau,N_\tau)].
\end{align}
We solve the stopping problem \eqref{opti pb} by a verification method: we first introduce the value function and determine the HJB equation. The HJB equation turns out to be an infinite system of ODEs for which it is difficult to find an explicit solution. We show via a fixed point argument that the ODE has a bounded solution on a small street interval. Via a continuation method, similar to the one used in \cite{rus2019some},  we then extend the solution to the whole street interval $[S, 0]$.  Finally, following a classical verification procedure, we identify the solution of the HJB equation with the value function. Thus, one can characterize an optimal stopping rule in terms of the solution function of the HJB equation. 
This optimal stopping time makes immediate use of new incoming information on the unknown intensity $\Lambda$. In this sense, the optimal stopping time can be interpreted as an online statistical learning decision policy.

\textbf{Comparison with the literature. } 
There is an extensive literature on the parking problem formulated as a stopping problem, see e.g.\ \cite{MacQueen} for an early reference in discrete time. A version of the parking problem in continuous time has been considered in \cite{sakaguchi1982optimal}. We refer to chapter 2.5 in \cite{ferguson2018optimal} for a recent overview on the variants of the parking problem that have been solved already.

\cite{Tamaki_1985} studies a Bayesian model for the parking problem in \emph{discrete time}. Under the assumption that the costs are given by the distance to some target, Tamaki falls back on discrete time dynamic programming to characterize optimal stopping rules. In particular, Tamaki observes that the One-Step-Look-Ahead rule is not optimal. The continuous-time framework of the present article opens the use of additional analytic methods for characterizing optimal stopping rules and hence facilitates to consider a larger class of cost functions. In particular, we allow for cost functions that depend on the number of observed free lots. Such cost functions are needed, e.g., when a group of drivers is searching several free parking lots so as to minimize their total costs (see Example \ref{ex2} below). 

If the parking problem needs to be solved in many consecutive rounds (e.g., every day), then one can assess a stopping and learning policy with the regret, the expected difference of the costs accumulated by the applied policy and the costs under the optimal stopping rule. \cite{ankirchner2026learning} propose a specific policy, based on a frequentist estimation of the possibly non-homogeneous jump intensity, for stopping round by round. They show that the policy has a regret that grows logarithmically with the number of rounds, and that there is no policy with a regret growing slower than the logarithm. 

\section{Model and main results}\label{sec:Static case Description}
Recall that  \(f_{\alpha,\beta}(x)=\frac{\beta^\alpha}{\Gamma(\alpha)}x^{\alpha-1}e^{-\beta x}, x\in[0,\infty)\), is the density of the gamma distribution with shape parameter \(\alpha >0\) and rate parameter \(\beta >0\). 

Let \(S\in (-\infty, 0), \alpha_S \in (0,\infty)\) and \(\beta_S \in (0,\infty)\). 
Let $\Lambda$ be gamma distributed with parameters $\alpha_S$ and $\beta_S$, and let \(N:=(N_t)_{t \in [S,\infty)}, N_S=0,\) be a homogeneous Poisson process with random jump intensity $\Lambda$, defined on the same probability space $(\Omega, \cF, P)$. 
We interpret the jump times of $N$ as the arrivals of free parking lots. We assume that a driver can stop at jump times of $N$ only. By observing the Poisson process the driver continuously updates the conditional distribution of $\Lambda$. 
Recall that the gamma distribution is the conjugate prior of the Poisson distribution. Hence the posterior distribution of $\Lambda$ at $t \in (S, \infty)$, given the information $\cF^N_t := \sigma(N_r: r \in [S, t])$, is the gamma distribution with parameters 
\begin{align*}
    \alpha_t := \alpha_S + N_t \quad \text{ and } \quad \beta_t := \beta_S + t-S.
\end{align*}
We define the set of admissible stopping times as the set 
\begin{align*}
    \cT:=\left\{\tau: \Omega \to [S, \infty] \right| \left.\tau \text{ is an } (\cF^N_t)_{t \in [S, \infty)}-\text{stopping time s.t. } \Delta_\tau=1 \right\}, 
\end{align*}
where $\Delta_t:=N_t-N_{t-}$.
    
Let $c:[S, \infty)\times \IN_0 \to [0, \infty)$ be a measurable cost function such that 
    \begin{enumerate}[label=(C\arabic*)]
         \item  $c(t,\cdot)$ is non-increasing for all $t \in[0, \infty)$, \label{C1}
         \item  $\exists K \in (0, \infty)$ and $ \delta \in [0, \alpha_S)$ such that for all $t\in [S, \infty)$ and $k\in\IN_0$ we have $ c(t, k)\leq K (1+|t|^\delta)$, \label{C2}
         \item  $c(t,N_t)$ is an $(\cF^N_t)$-submartingale on $[0, \infty)$ such that for $\tau,\si\in\cT$ with $\tau\geq\sigma\geq0$ a.s. we have $\IE[c(\tau,N_\tau)]\geq \IE[c(\sigma,N_\sigma)]$. \label{C3}
     \end{enumerate}
    The car driver seeks to minimize the costs \(c(\tau,N_\tau)\) corresponding to the chosen parking lot \(\tau \in \cT\).
    Such cost can be independent of the second argument, e.g., $c(t,k)=|t|$ (see Example~\ref{example_cost_absolute} below)
    or can depend on $N_\tau$, e.g., $c(t,k)=|t|+\IE[\inf\{|s|:s>t,\Delta_s=1\}\mid N_t=k]$ (see Example~\ref{ex2} below). 
    We discuss the Assumptions \ref{C1}--\ref{C3} in more detail in Remarks~\ref{rem:Value t>0} and~\ref{rem:t>0 Stoppzeit<tau_t^1} below. 

In order to determine a stopping time minimizing the expected costs, we introduce the value function and characterize it in terms of the associated HJB equation. To this end, for all $t\in [S, \infty)$, we denote by $\cT_t = \{\tau \in \cT: \tau \ge t\}$ the set of admissible stopping times with values greater than or equal to $t$. For every $(t, k, \Delta) \in [S, \infty) \times \N_0 \times \{0, 1\}$ the expected costs of a stopping time $\tau\in \cT_t$ conditional on the event $\{N_t = k\} \cap \{\Delta_t = \Delta\}$ are given by
    \begin{align}\label{eq:defiJ}
        J(t,k,\Delta,\tau) := \IE[c(\tau,N_\tau)|N_t = k, \Delta_t=\Delta].
    \end{align}
    (Note that $\{\Delta_t=1\}$ is a $P$-null set. Hence, conditioning on $\{\Delta_t=0\}$ can be omitted in the definition~\eqref{eq:defiJ}.)
   The value function is defined by the minimal expected costs:
    \begin{align*}
        V(t,k,\Delta):=\inf_{\tau\in \cT_t} J(t,k,\Delta, \tau),\quad \quad (t,k,\Delta)\in [S,\infty)\times \IN_0\times \{0,1\}.
    \end{align*}
    So the driver's state has three dimensions: $t\geq S$ is the current position on the street, $k\in\IN_0:=\{0,1,2,...\}$ is the number of observed free parking lots until time $t$, and $\Delta\in\{0,1\}$ indicates the state of the current parking lot: 1 means free, 0 means  occupied.
    \begin{rem}\label{rem:Value t>0}
         From Condition~\ref{C3} follows that it is optimal to stop as early as possible after zero, provided one has not stopped before.
    Therefore it is enough to consider stopping times $\tau\in \cT_t: \tau\leq\max\{\tau_0^1,\tau_t^1\}$, where $\tau_t^i:=\inf\{s> t:N_s-N_t=i\}$ is the $i$-th jump time of $N$ after time $t$. 
    Furthermore, the value function satisfies
    \begin{align}\label{eq:Boundary condi}
     V(t,k,1)&=c(t,k), && \forall t\geq 0, k \in \N_0,\\\label{eq:Boundary condii}
     V(t,k,0)&=\IE  [c(\tau_t^1,N_{\tau_t^1})|N_t = k],&&\forall t\geq 0, k \in \N_0.
    \end{align}
    \end{rem}
We will characterize the value function as the unique bounded solution to the associated HJB equation. A crucial ingredient of the HJB equation is the compensator of $N$ with respect to the observer's filtration $(\cF^N_t)$, explicitly provided in the next proposition. 
\begin{prop}
The compensator of $N$ with respect to the filtration $(\cF^N_t)_{t \in [S, \infty)}$ is given by 
\begin{align*}
\int_S^t \lambda(u,N_{u-}) du, \quad t \in [S, \infty),
\end{align*}
where \begin{align}
    \lambda(t,k):= \frac{\alpha_S + k}{\beta_t}, \text{ for all } (t,k) \in [S,\infty) \times \N_0.  
\end{align}
\end{prop}
\begin{proof}
We need to show that for $u, t \in [S, \infty)$ with $u < t$ we have 
\begin{align}\label{mart equ}
    \IE\left[ N_t - N_u | \cF^N_u \right] = \IE\left[ \int_u^t  \lambda(r,N_{r-}) dr\middle| \cF^N_u \right].  
\end{align}
The LHS of \eqref{mart equ} satisfies 
\begin{align*}
 \IE\left[ N_t - N_u \middle| \cF^N_u \right] & =  (t-u) \IE(\Lambda | \cF^N_u) = (t-u) \int_0^\infty x f_{\alpha_{u}, \beta_u} (x) dx \\
 & = (t-u)  \frac{\alpha_{u}}{\beta_u}= (t-u)\lambda(u, N_u).     
\end{align*}
For the RHS of \eqref{mart equ} we have 
\begin{align*}
    \IE\left[ \int_u^t  \lambda(r,N_{r-}) dr \middle| \cF^N_u \right] & = \IE\left[ \int_u^t  \frac{\alpha_{r-}}{\beta_r} dr \middle| \cF^N_u \right] = \IE\left[ \int_u^t  \frac{\alpha_r}{\beta_r} dr \middle| \cF^N_u \right]\\
    & = \int_u^t  \IE\left[ \frac{\alpha_r}{\beta_r} \middle| \cF^N_u \right] dr = \int_u^t   \frac{\alpha_S + \IE\left[N_{r} \middle| \cF^N_u \right]}{\beta_r} dr. 
\end{align*}
Since $\IE\left[N_{r}| \cF^N_u \right] = N_u + (r-u) \IE(\Lambda | \cF^N_u) = N_u + (r-u)\lambda(u, N_u)$ we further get
\begin{align*}
    \IE\left[ \int_u^t  \lambda(r,N_{r-}) dr \middle| \cF^N_u \right] & = \int_u^t \frac{\alpha_S + N_u + (r-u)\lambda(u, N_u)}{\beta_r} dr \\
    & = \int_u^t \frac{\alpha_S + N_u }{\beta_r}  + \frac{(r-u)(\alpha_S+ N_u)}{\beta_r \beta_u} dr \\
     & = \int_u^t \frac{\alpha_S + N_u }{\beta_r}  + 
     (\alpha_S+ N_u)\left(\frac{1}{\beta_u}-\frac{1}{\beta_r}\right) dr \\
    & = \int_u^t \frac{\alpha_S+N_u}{\beta_u}   dr \\
    & = (t-u) \lambda(u, N_u). 
\end{align*}
Thus we have shown \eqref{mart equ}. 
\end{proof}
The HJB equation uniquely characterizing the value function of our stopping problem is the following system of equations, indexed by $k \in \N_0$,
\begin{align}
\begin{split}\label{eq:HJB} 
        \partial_t W(s,k,0)&=\la(s,k) (W(s,k,0)-W(s,k + 1,1)), 
        \\
        W(s,k,1)&=
        \min\{c(s,k), W(s,k,0)\}, \quad s\in [S,0].
\end{split}
    \end{align}
Our main theoretical result is the following.  
    
\begin{thm}\label{thm:valueverifi}
    The value function $V$ is the unique bounded solution of \eqref{eq:HJB} on $[S,0]$ for all $k \in \N_0$ subject to the boundary conditions \eqref{eq:Boundary condi} and \eqref{eq:Boundary condii}.
    An optimal stopping time is given by \begin{align}\label{optiStop}
        \tau^*:=\inf \{t\geq S: V(t,N_t, 0)\geq c(t,N_t), \Delta_t=1\}.
    \end{align}
\end{thm}
To prove the theorem, we first show existence of a bounded solution of the HJB equation, see Section~\ref{sec:existence} below, and then verify that the solution coincides with the value function, see Section~\ref{sec:verification}.

We close this section with some remarks on the assumptions on the cost function and with several examples of cost functions satisfying the assumptions. To shorten notation, we denote by $P_{t,k} = P (\cdot| N_t = k)$ the probability measure $P$ conditioned on the event $\{N_t = k\}$, for all $(t,k) \in [S, \infty) \times \N_0$. Moreover, we use $\tau_t^i:=\inf\{s> t:N_s-N_t=i\}$, cf.\ Remark~\ref{rem:Value t>0}.
\begin{rem}\label{rem:t>0 Stoppzeit<tau_t^1}
    Condition \ref{C1} assures that observing more parking lots does not lead to higher expected costs.
    Condition \ref{C2} guarantees that there exists a stopping time with finite expected costs (see Lemma~\ref{lem:Bounded boundary} below).     
    Conditions \ref{C1} and \ref{C3} together imply monotonicity of $c(\cdot,k)$ on $[0,\infty)$:
    \begin{align}\label{eq:monotony in t}
        c(t,k)\leq \IE_{t,k}[c(t+h,N_{t+h})]\leq\IE_{t,k}[c(t+h,k)]= c(t+h,k).
    \end{align}
    \end{rem}
\begin{exa}
\label{example_cost_absolute}
Let $c(t,k) = |t|$ for all $k \in \N_0$. Note that Condition \ref{C1} is trivially satisfied.  Condition \ref{C2} is only satisfied if $\alpha_S > 1$. Notice that  $\alpha_S > 1$ guarantees that the expected costs $\IE_{0,k}|\tau_0^1|$ are finite for all $k$ (cf.\ Lemma~\ref{lem:Bounded boundary}). Finally, it is immediate to see that $c(t,N_t)=|t|$ is a submartingale and monotonously increasing on $[0,\infty)$ and hence also Condition \ref{C3} is satisfied. 
\end{exa}
\begin{exa} \label{ex2} 
Let $c(t,k)=|t|+\IE_{t,k}|\tau_t^1|$, the expected total distance to the target location of two consecutive free parking lots when the first free lot is at $t$. 
Condition \ref{C1} is entailed by Lemma \ref{lem:monotony in k of E} below. 
Condition \ref{C2} is satisfied if $\alpha_S > 1$ since $\tau_t^1-t$ is Lomax distributed  (cf. Lemma~\ref{lem:Lomax}) whose mean  $\IE_{t,k}(\tau_t^1-t)$ is bounded uniformly in $k$. 
Hence there exists a $K$ such that $c(t,k)\leq 2|t|+\IE_{t,k}|\tau_t^1-t|\leq K(1+|t|)$.
To show that $c(t,N_t)$ is a submartingale on $[0,\infty)$ it is enough to show that
        \begin{align*}
            0\leq \partial_t c(t,k) +\la(t,k) (c(t,k+1)-c(t,k)).
        \end{align*}
    Indeed, by Lemma~\ref{lem:Lomax} below and the properties of the Lomax distribution, we have
        \begin{align*}
            0&\leq 
            \partial_t c(t,k) +\la(t,k) (|t|+\IE_{t,k+1}|\tau_t^1|-|t|-\IE_{t,k}|\tau_t^1|)\\
            &= 
            \partial_t \left(t+t+\frac{\beta_t}{\alpha_S+k-1}\right) 
                +\frac{\alpha_S+k}{\beta_t} 
                \left(\frac{\beta_t}{\alpha_S+k}-\frac{\beta_t}{\alpha_S+k-1}\right)\\
            &=2+\frac{1}{\alpha_S+k-1} +1-\frac{\alpha_S+k}{\alpha_S+k-1}\\
            &=2.
        \end{align*}
     Condition~\ref{C3} holds, since beside the submartingale property, we have also the following for all $\tau\geq\si\ge0$ a.s.:
     \begin{align}
         \IE[c(\tau,N_\tau)]=\IE\tau + \IE\tau_\tau^1\geq\IE\si+\IE\tau_\si^1=\IE [c(\si,N_\si)].
     \end{align}
    To have finite expected costs $\IE_{0,k}[|\tau_0^1|+|\tau_0^2|]$ for all $k$, a prior parameter $\alpha_S>1$ is required.
    \end{exa}
    \begin{exa}
    Let $c(t,k) = |t| + \sI_A (t), A\subset [S,0]$. The indicator set $A$ can be interpreted as an area where parking charges apply. Note that Conditions \ref{C1}--\ref{C3} are satisfied if $\alpha_S>1$ (analogously to Example~\ref{example_cost_absolute}). 
    \end{exa}
\section{
Existence of a solution of the  HJB equation}\label{sec:existence}
In this section we prove that there exists a bounded solution of the HJB equation (see Equations~\eqref{eq:HJB}). 
First, we show the following for a truncated system.
\begin{lemm}\label{lem:uniquePDE_truncated}
    There exists a bounded, non-negative function $w: [S,0]\times \IN_0\to \IR$ that is absolutely continuous in its first variable, satisfies for a.e.\ $ s \in [S,0]$ 
    \begin{align}\label{hjb_truncated}
        \partial_t w(s,k)&=\la(s,k) (w(s,k)-\min\{c(s,k+1), w(s,k+1)^+\}),\quad \forall k\in \N_0,
    \end{align}
and satisfies the boundary condition  
\begin{align}\label{boundary condi_truncated}
    w(0,k)=\IE_{0,k}  [c(\tau_0^1,k+1)] \text{ for all } k \in \N_0.
\end{align} 
\end{lemm}
For proving Lemma \ref{lem:uniquePDE_truncated} we need several auxiliary results.
\begin{lemm}\label{lem:Lomax}
   The distribution  of \(\tau_t^1-t\) conditional on $\cF_t^N$ is $\mathrm{Lomax}(\alpha_t,\beta_t)$.
   More precisely, the conditional distribution  of \(\tau_t^1-t\) has a density
   given by
   \begin{align*}
     p_{\tau_t^1-t | \cF_t^N}(s)=\frac{\alpha_t \beta_t^{\alpha_t}}{(\beta_t+s)^{\alpha_t+1}}, \quad s\ge0.  
   \end{align*}
\end{lemm}
\begin{proof}
It is a standard result that the expected stochastic intensity $\IE_{t,k}[\lambda(s,N_{s-})|\tau_t^1-t>s],s>0,$ coincides with the hazard rate $h_t(s):=\lim_{\Delta s \to 0^+} \frac{P(\tau_t^1 -t\le s+\Delta s|\tau_t^1-t >s)}{\Delta s}$ of the waiting time until the next jump  (see, e.g.,  \cite[theorem T7]{brémaud1981point}). 
Together with the identity $P_{t,k}(\tau_t^1-t>s) 
    =\exp\left(-\int_t^{t+s} h_t(u) \td u\right)$, a classical result from survival analysis, we have 
    \begin{align}\label{surv prob}
    P_{t,k}(\tau_t^1-t>s) 
    &=\exp\left(-\int_t^{t+s} \IE_{t,k}(\la(u,N_{u-})|\tau_t^1>u) \td u\right) \nonumber \\
    &=\exp\left({-\int_t^{t+s} \frac{\alpha_S+k}{\beta_u} \td u }\right) \nonumber \\
    &=\exp(-(\alpha_S+k) [\ln(\beta_u)]_{u=t}^{u=t+s})
    \nonumber \\
    &
    = \left(\frac{\beta_t}{\beta_{t+s}}\right)^{\alpha_S+k}.
\end{align}
By taking the derivative we obtain the density of \(\tau_t^1-t\):
\begin{align*}
    p_{\tau_t^1-t | \cF_t^N}(s)
    &=- \frac{\td}{\td s}P(\tau_t^1-t>s|\cF_t^N) 
    = \frac{\alpha_t \beta_t^{\alpha_t}}
            {\beta_{t+s}^{\alpha_t+1}}.
\end{align*}
\end{proof}
\begin{lemm}\label{lem:monotony in k of E}
    The sequence $\{\IE_{0,k}  [c(\tau_0^1,k+1)]\}_{k\in\IN_0}$ is non-increasing in $k$.
\end{lemm}
\begin{proof}
    Let $k\in\IN_0$ and $F_k(t)$ be the distribution of the first jump time of $N$ after zero, conditional on $N_0=k$.
    Then $F_k(t)\leq F_{k+1}(t)$ for all $t\geq0$, because the jump intensity $\la$ on $[0,\infty)$, conditional on $N_0=k+1$, is larger than the jump intensity conditional on $N_0=k$.
    Therefore, since $c(\cdot,k)$ is non-decreasing on $[0,\infty)$, 
    \begin{align*}
        \int_0^\infty c(t,k+1)\td F_k(t)
        \geq \int_0^\infty c(t,k+1)\td F_{k+1}(t)
        \stackrel{\tiny \ref{C1}}\geq \int_0^\infty c(t,k+2)\td F_{k+1}(t).
    \end{align*}
\end{proof}
\begin{lemm}\label{lem:Bounded boundary}
    The sequence $\{\IE_{0,k}  [c(\tau_0^1,k+1)]\}_{k\in\IN_0}$ is bounded.
\end{lemm}
\begin{proof}
By Lemma~\ref{lem:monotony in k of E} we have $\IE_{0,k}  (c(\tau_0^1,k+1))\leq \IE_{0,0}  (c(\tau_0^1,1))$. Moreover, using the abbreviations $\alpha = \alpha_S$ and $\beta = \beta_0 = \beta_S-S$, we obtain with Condition \ref{C2}
\begin{align*}
    \IE_{0,0} (c(\tau_0^1,1))&=\int_0^\infty c(t,1) \frac{\alpha \beta^\alpha}{(\beta+t)^{\alpha+1}} \,\td t
\\
&\leq
    K\alpha\beta^\alpha\int_0^\infty (1+t^\delta ) \frac{1}{(\beta+t)^{\alpha+1}} \, \td t\\
    &\leq K\alpha\beta^\alpha\int_0^\infty (\beta+t)^{-(\alpha+1)} +
    (\beta+t)^{\delta-\alpha-1} \,\td t<\infty.
\end{align*}
Note that the integral is finite, because $\delta-\alpha < 0$.
Hence we have shown the claim. 
\end{proof}
\begin{lemm}\label{lemma existence}
Let $a \in (S, 0]$. Suppose that there exists $w^*:[a,0] \times \N_0 \to \R$ that solves \eqref{hjb_truncated} on $[a,0]$ and satisfies $w^*(0,k) = \IE_{0,k}  [c(\tau_0^1,k+1)]$ and $\sup_{k \ge 0} |w^*(a,k)| < \infty$.
Then there exists $b < a$ such that one can extend $w^*$ to a solution of \eqref{hjb_truncated} on the interval $[b,0]$. Moreover, there exists $h > 0$, independent of $a$, such that one can choose $b$ with $a-b \ge \min\{ h, a-S\}$.      
\end{lemm}
\begin{proof}
Let $w^*:[a,0] \times \N_0 \to \R$ be an absolutely continuous solution of \eqref{hjb_truncated} on $[a,0]$ satisfying $w^*(0,k) = \IE_{0,k}  [c(\tau_0^1,k+1)]$ and $\sup_k |w^*(a,k)| < \infty$. For $\rho >0$ and $b \in [S, a)$ and all functions $v:[b,a] \times \N_0 \to \mathbb{R}$, we define 
\[||v||_\rho^{b,a}:=\sup_{t \in [b,a]} \sup_{k\in \mathbb{N}_0} e^{-\rho (\alpha_S +k)g_a(t)} |v(t,k)|, 
\]
where \(g_a(t):=
\ln (\beta_a/\beta_t).\) 
Note that the set
$$B_\rho(b,a):=\{v: [b,a] \times \mathbb{N}_0 \to \mathbb{R} \mid v \text{ is measurable and } ||v||_\rho^{b,a} < \infty\} $$   
is a Banach space. For any $v \in B_\rho(b,a)$ we define a new function $\phi(v)$ with domain $[b,a] \times \N_0$ by 
\begin{align}\label{eq:defOperatorPhi}
    \phi(v)(t,k):=
     w^*(a,k)  +\int_t^a \lambda(s,k)[\min\{ c(s,k+1),v(s,k+1)^+\}-v(s,k)]\td s.
\end{align}
Note that \(|\min\{z,x^+\}-\min\{z,y^+\}|\le |x^+-y^+|\le |x-y|\) for \(x,y,z \in \mathbb{R}\). For $v, w \in B_\rho(b,a)$ we have 
\begin{align*}
        |\phi(v)(t,k)-\phi(w)(t,k)| &\leq \int_t^a \lambda(s,k) \big|\left(\min\{c(s,k+1), v(s,k+1)^+\}-v(s,k)\right) \\
    &\hspace{3cm} -\left(\min\{c(s,k+1), w(s,k+1)^+\}-w(s,k)\right)\big| \, \td s\\
    &\le \int_t^a \lambda(s,k) (| v(s,k+1)-w(s,k+1) | + |v(s,k)-w(s,k)|) \, \td s \\
    &= \int_t^a \lambda(s,k) e^{-\rho (\alpha_S +k)g_a(s)}e^{\rho (\alpha_S +k)g_a(s)}(| v(s,k+1)-w(s,k+1) | \\& \hspace{6.5cm} + |v(s,k)-w(s,k)|) \, \td s \\ 
    &\le 2 ||v-w||_\rho^{b,a} \int_t^a e^{\rho (\alpha_S +k+1) g_a(s)} \lambda(s,k) \, \td s \\
    &\le 2 ||v-w||_\rho^{b,a} \int_t^a e^{\rho (\alpha_S +k+1) g_a(s)}(\alpha_S +k+1)(-g_a'(s)) \, \td s. 
    \end{align*}
    By substituting \(u=g_a(s)\)
    we get 
    \begin{align*}
        |\phi(v)(t,k)-\phi(w)(t,k)|&\leq -2 ||v-w||_\rho^{b,a} \int_{g_a(t)}^{g_a(a)} e^{\rho (\alpha_S +k+1) u}(\alpha_S +k+1) \, du \\
&=2 ||v-w||_\rho^{b,a} \left[\frac1\rho e^{\rho (\alpha_S +k+1) u}\right]_0^{g_a(t)}  \\
        &=2 ||v-w||_\rho^{b,a} \frac1\rho \left[e^{\rho(\alpha_S +k+1) g_a(t)}-1\right], 
    \intertext{and hence }
         e^{-\rho (\alpha_S  + k)g_a(t)}
         |\phi(v)(t,k)-\phi(w)(t,k)|&\leq 2 ||v-w||_\rho^{b,a} \frac1\rho \left[e^{\rho g_a(t)}-e^{-\rho(\alpha_S +k)g_a(t)}\right]\\  
         &\le \frac2\rho e^{\rho g_a(t)} ||v-w||_\rho^{b,a}.
    \end{align*}
Since the last inequality holds for arbitrary \((t,k)\) we can conclude
\begin{align}
\label{contraction}
    ||\phi(v)-\phi(w)||_\rho^{b,a} &\le  \sup_{t\in [b,a]}\frac2\rho e^{\rho g_a(t)} ||v-w||_\rho^{b,a}=\frac2\rho e^{\rho g_a(b)}||v-w||_\rho^{b,a}.
\end{align}
We use inequality \eqref{contraction} first to show that $\phi$ is an operator from $B_\rho(b,a)$ to itself. Indeed, 
\begin{align*}
    \|\phi(v)\|_\rho^{b,a} \le \|\phi(v) - \phi(0)\|_\rho^{b,a} + \|\phi(0)\|_\rho^{b,a} \le \frac2\rho e^{\rho g_a(b)}||v||_\rho^{b,a} + \sup_k |w^*(a,k)| < \infty. 
\end{align*}
Now suppose that $\rho > 2$. Then we can choose $b \in [S, a)$ such that $\frac2\rho e^{\rho g_a(b)} < 1$. Inequality~\eqref{contraction} implies that $\phi$ is a contraction on $B_\rho(b,a)$ in this case. Hence, by the Banach fixed point theorem, there exists a unique $v^* \in B_\rho(b,a)$ such that $\phi(v^*)(t,k) = v^*(t,k)$ for all $(t,k)$. Now extend $w^*$ to $[b,0]$ by setting $w^*(t,k) = v^*(t,k)$ for all $t \in [b,a)$ and $k \in\N_0$. 

Notice that $w^*$ is a solution of \eqref{hjb_truncated} on the extended interval $[b,0]$. Indeed, by the very construction we have for all $t \in [b,0]$ and $k \in \N_0$
\begin{align*}
w^*(t,k) = \IE_{0,k}  [c(\tau_0^1,k+1)]  +\int_t^0 \lambda(s,k)[\min\{w^*(s,k+1)^+, c(s,k+1)\}-w^*(s,k)]\td s. 
\end{align*}
This shows that $w^*$ is absolutely continuous in $t$, satisfies the ODE \eqref{hjb_truncated} for a.e. \(t\) and the terminal condition $w^*(0,k) = \IE_{0,k}  [c(\tau_0^1,k+1)]$. 

Finally, we show the last statement of the lemma. To this end we seek an $h>0$ such that for all $a \in [S+h, 0]$
\begin{align}\label{condi h}
\frac{2}{\rho}e^{\rho g_{a}(a-h)} = \frac2\rho \exp\left(\rho \ln\left(\frac{\beta_a}{\beta_{a-h}}\right)\right) \overset{!}{<}1. 
\end{align}
Note that \eqref{condi h} is fulfilled if $\rho > 2$ and \(h < h^*: = \beta_S (1- (\frac{\rho}{2})^{-1/\rho})\). 
\end{proof}
\begin{rem}
    The statement of the previous lemma  holds true for $a=0$. In this case only the boundary condition~\eqref{boundary condi_truncated} has to be satisfied.
\end{rem}
\begin{proof}[Proof of Lemma \ref{lem:uniquePDE_truncated}]
We apply Lemma~\ref{lemma existence} successively from zero to $S$. 
Let initially $a_0=0$ and $b_0=-h$. In each step $i$, we extend the interval by a length of $h$ by choosing its boundaries $a_i=a_{i-1}-h, b_i=\max\{S, a_i-h\}$, $i \in \{1,2,..., \lceil |S|/h\rceil-1\}$.
To this end, it remains to show that each extension is bounded.
We do this by contradiction in the next paragraph.
The last extension is then
a bounded function $w: [S, 0] \times \N_0 \to \R$ satisfying the ODE \eqref{hjb_truncated} and the boundary condition $w(0,k)=\IE_{0,k}  [c(\tau^1_0,k+1)]$.  

From Lemma~\ref{lem:Bounded boundary} we know that $\sup_k \IE_{0,k} (c(\tau_0^1,k+1))$ is bounded, say by $\bar K$. Let \(M:= \max\{K(1+|S|^\delta),\bar K\}.\) By (C2) we have \(c(t,k)\le K(1+|S|^\delta)\le M\) for all \((t,k)\in [S,0]\times \mathbb{N}_0.\)
\begin{enumerate}
    \item \textit{Upper Bound:}
    We assume that there exists \(k_0 \in \mathbb{N}_0\) and \(t_0 \in [b_i,0]\) such that \(w(t_0,k_0) > M.\)  Let \(t^*:=\sup\{t \in [b_i,0]: w(t,k_0)\ge \frac{w(t_0,k_0)+M}{2} \}.\) We know that \(t^*<0\) because \(\frac{w(t_0,k_0)+M}{2}>M\ge w(0,k_0)=\IE_{0,k_0}[c(\tau_0^1,k_0+1)]\).  Since \(w(\cdot, k_0)\) is continuous it holds \(w(t^*,k_0)=\frac{w(t_0,k_0)+M}{2}\) and \(w(t,k_0) < \frac{w(t_0,k_0)+M}{2}\) for \(t\in (t^*,0].\) Again by continuity of \(w(\cdot, k_0)\) there exists an \(\epsilon >0\) such that \(w(\hat t, k_0) >M\) for all \(\hat t \in [t^*,t^*+\epsilon]\). It follows for a.e. \(\hat t \in [t^*,t^*+\epsilon]\) that \begin{align*}
        \partial_t w(\hat t , k_0)&=\lambda (\hat t ,k_0)(w(\hat  t,k_0)-\min\{c(\hat t, k_0+1), w(\hat t, k_0+1)^+\})\\
        &>0,
    \end{align*}
    where the latter inequality follows from \(w(\hat  t,k_0)> M\) and \(\min\{c(\hat t, k_0+1), w(\hat t, k_0+1)^+\} \le M.\) 
    Since \(w(\cdot, k_0)\) is absolutely continuous, we obtain 
    \begin{align*}
       w(\hat t, k_0)=w(t^*, k_0)+ \int_{t^*}^{\hat t} \partial_s w(s,k_0) \, \td s > w(t^*, k_0)= \frac{w(t_0,k_0)+M}{2}
    \end{align*}
    for all \(\hat t \in (t^*,t^*+\epsilon]\).
   From the definition of $t^*$, however, it must hold $w(\hat t, k_0) < \frac{w(t_0,k_0)+M}{2}$, and hence we obtain a contradiction. 

    \item \textit{Lower Bound:}
    We assume that there exists \(k_0 \in \mathbb{N}_0\) and \(t_0 \in [b_i,0]\) such that \(w(t_0,k_0) <0.\)  Let \(t^*:=\sup\{t \in [b_i,0]: w(t,k_0)\le \frac{w(t_0,k_0)}{2} \}.\) We know that \(t^*<0\) because \(\frac{w(t_0,k_0)}{2}<0\le w(0,k_0)=\IE_{0,k_0}[c(\tau_0^1,k_0+1)]\).  Since \(w(\cdot, k_0)\) is continuous it holds \(w(t^*,k_0)=\frac{w(t_0,k_0)}{2}\) and \(w(t,k_0) > \frac{w(t_0,k_0)}{2}\) for \(t\in (t^*,0]\).
    Again by continuity of \(w(\cdot, k_0)\) there exists an \(\epsilon >0\) such that \(w(\hat t, k_0) <0\) for \(\hat t \in [t^*,t^*+\epsilon]\).  It follows for a.e. \(\hat t \in [t^*,t^*+\epsilon]\) that \begin{align*}
        \partial_t w(\hat t , k_0)&=\lambda (\hat t ,k_0)(w(\hat  t,k_0)-\min\{c(\hat t, k_0+1), w(\hat t, k_0+1)^+\})\\
        &<0,
    \end{align*}
    where the latter inequality follows from \(w(\hat  t,k_0)<0\) and \(\min\{c(\hat t, k_0+1), w(\hat t, k_0+1)^+\} \ge0.\) Since \(w(\cdot, k_0)\) is absolutely continuous, we obtain  \(w(\hat t, k_0)<\frac{w(t_0,k_0)}{2}\) for all \(\hat t \in (t^*,t^*+\epsilon]\), which leads to a contradiction. 
\end{enumerate}
Since \(k_0\) was arbitrary, we have \(0 \le w \le M.\)
\end{proof}

\begin{thm}\label{lem:uniquePDE}
    There exists a bounded function $w: [S,0]\times \IN_0\to \IR$ that is absolutely continuous in its first variable, satisfies for almost every $s \in[S,0]$ 
    \begin{align}\label{hjb for 0 comp}
        \partial_t w(s,k)&=\la(s,k) (w(s,k)-\min\{c(s,k+1), w(s,k+1)\}),  \text{ for all }k\in \N_0,
    \end{align}
and satisfies the boundary condition  
\begin{align}\label{boundary condi}
    w(0,k)=\IE_{0,k}  [c(\tau_0^1,k+1)] \text{ for all } k \in \N_0.
\end{align} 
\end{thm}

\begin{proof}
   Let  \(w\) be the function from Lemma \ref{lem:uniquePDE_truncated} which solves \eqref{hjb_truncated}.  Since \(w\ge 0\), we have \(w^+=w\). Therefore, the right-hand sides of \eqref{hjb_truncated} and \eqref{hjb for 0 comp} coincide for this function \(w\). Moreover, \(w\) is bounded, absolutely continuous in its first variable and satisfies \eqref{boundary condi}. Hence \(w\) is a solution of \eqref{hjb for 0 comp} with the asserted properties.
\end{proof}

\section{Verification}\label{sec:verification}
In this section, we verify that any bounded solution of the HJB equation \eqref{eq:HJB} satisfying the boundary conditions \eqref{eq:Boundary condi} and \eqref{eq:Boundary condii} coincides with the value function. In particular, this shows that the value function can be characterized as the \emph{unique} bounded solution of the HJB equation. 

For the verification we use the following lemmata. We use the notation $s\wedge t:=\min\{s,t\}$.

\begin{lemm} \label{lem:sum1}
Let $W: [S,0]\times \IN_0\times \{0,1\}\to \IR$ be a bounded function and $\si\in \cT_t^0:=\{\tau\wedge0:\tau\in \cT_t\}$.
Let $k^*:=\min\{j: \tau_t^j\wedge0=\si\}$ be the number of jumps between $t$ and the stopping time $\si$ (zero is counted as jump in the definition of $k^*$).
Moreover let $\sigma_t^0:=t$ and $\sigma_t^j:=\tau_t^j\wedge 0$ be the $j$-th jump time of $N$ after $t$ truncated at zero.
Then
\begin{align}
    W(\si,N_\si,\Delta_{\sigma})&= 
    W(t,N_t,0)
        +\sum_{1\leq j\leq k^*} \big\{ W(\sigma_t^j,N_{\sigma_t^j},\Delta_{\sigma_t^j})-W(\sigma_t^j,N_{\sigma_t^{j-1}},0) \big\}    \label{eq:line1}\\ &\quad 
        +\sum_{1\leq j\leq k^*} \big\{ W(\sigma_t^j,N_{\sigma_t^{j-1}},0)-W(\sigma_t^{j-1},N_{\sigma_t^{j-1}},0) \big\}
    \label{eq:line2}\\ &\quad +
    \sum_{2\leq j\leq k^*} \big\{ W(\sigma_t^{j-1},N_{\sigma_t^{j-1}},0)-W(\sigma_t^{j-1},N_{\sigma_t^{j-1}},1) \big\},\quad a.s.\label{eq:line3}
\end{align}
\end{lemm}
\begin{proof}
    Since the sums are finite almost surely, they are interchangeable. Note that
    \begin{align*}
    W(\si,N_\si,\Delta_\si)
    &=
     W(t,N_t,0)+ W(\si,N_{\si},\Delta_\si)
    -W(t,N_t,0)
    \\
    &=
     W(t,N_t,0)+ W(\si,N_{\si},\Delta_\si)
    \\ &\quad +
    \sum_{1\leq j\leq k^*} \big\{-W(\sigma_t^{j-1},N_{\sigma_t^{j-1}},0) \big\}
    \\ &\quad +
    \sum_{2\leq j\leq k^*} \big\{ W(\sigma_t^{j-1},N_{\sigma_t^{j-1}},0) \big\}\\
    &=    
    W(t,N_t,0)+\sum_{1\leq j\leq k^*} \big\{ 
        W(\sigma_t^j,N_{\sigma_t^j},\Delta_{\sigma_t^j})-{W(\sigma_t^j,N_{\sigma_t^{j-1}},0)} 
    \big\}
    \\ &\quad +
    \sum_{1\leq j\leq k^*} \big\{{ 
        W(\sigma_t^j,N_{\sigma_t^{j-1}},0)}-W(\sigma_t^{j-1},N_{\sigma_t^{j-1}},0) \big\}
    \\ &\quad +
    \sum_{2\leq j\leq k^*} \big\{ 
        W(\sigma_t^{j-1},N_{\sigma_t^{j-1}},0)-W(\sigma_t^{j-1},N_{\sigma_t^{j-1}},\Delta_{\sigma_t^{j-1}}) \big\} ,\quad a.s.
\end{align*}
Since $\Delta_{\sigma_t^{j}}=1$ for $j<k^*$, the result is shown.
\end{proof}
\begin{lemm}\label{lem:EW}
Let $W(\cdot,k,0)$ be bounded and absolutely continuous.
Let $k^*, \si$ and $\si_t^j$ be defined as in Lemma~\ref{lem:sum1}. Then
\begin{align}
    \IE_{t,k}  W(\si,N_\si,\Delta_{\sigma})&= W(t,k,0)+\IE_{t,k}  \int_t^\si \lambda(s,N_{s-})\big(W(s,N_{s-}+1,1)-W(s,N_{s-},0) \big) \td s
    \label{eq:line1equiv}\\ &\quad +
    \IE_{t,k}  \int_t^\si \partial_t W(s,N_s,0) \td s
    \label{eq:line2equiv}\\ &\quad +
    \IE_{t,k} 
    \sum_{2\leq j\leq k^*} \big\{ W(\sigma_t^{j-1},N_{\sigma_t^{j-1}},0)-W(\sigma_t^{j-1},N_{\sigma_t^{j-1}},1) \big\}.\notag
\end{align}
\end{lemm}
\begin{proof} 
We compare Lines \eqref{eq:line1} and \eqref{eq:line2} of Lemma~\ref{lem:sum1} with the terms in Lines~\eqref{eq:line1equiv} and \eqref{eq:line2equiv}, respectively.
Since $W $ is absolutely continuous in the first argument, we have equality of Lines~\eqref{eq:line2} and \eqref{eq:line2equiv}:
    \begin{align*}
        \int_t^\si \partial_t W(s,N_s,0) \td s &=
        \sum_{1\leq j\leq k^*}
        \int_{\sigma_t^{j-1}}^{\sigma_t^j} \partial_t W(s,N_s,0) \td s\\
        &=\sum_{1\leq j\leq k^*} \big\{ W(\sigma_t^j,N_{\sigma_t^{j-1}},0)-W(\sigma_t^{j-1},N_{\sigma_t^{j-1}},0) \big\}.
    \end{align*}
Further we get via compensation the equality of the terms in Lines~\eqref{eq:line1} and \eqref{eq:line1equiv}. 
For that, note that the integral of a bounded process with respect to a martingale is again a martingale.
Hence, $ \int_t^\si (W(s,N_{s-}+1,1)-W(s,N_{s-},0))(\td N_s-\lambda(s,N_{s-})\td s)$ is a martingale and
        \begin{align*}
            0&=\IE_{t,k}  \int_t^{\si} (W(s,N_{s-}+1,1)-W(s,N_{s-},0))(\td N_s-\lambda(s,N_{s-})\td s)
            \end{align*}
            by Doob's optional stopping theorem with bounded stopping times. 
 By the additivity of the integrator, we obtain 
            \begin{align}
            0=&\IE_{t,k}  \left(\sum_{1\leq j\leq k^*} \big\{ 
                W(\sigma_t^j,N_{\sigma_t^j},\Delta_{\sigma_t^j})-W(\sigma_t^j,N_{\sigma_t^{j-1}},0) \big\}\right)\notag\\
            &-\IE_{t,k}  \int_t^\si \lambda(s,N_{s-}) (W(s,N_s+1,1)-W(s,N_s,0))\td s.\label{eq:EMart=0}
        \end{align}
\end{proof}
\begin{rem}
In Lemma~\ref{lem:EW} the stopping time $\si$ is bounded from above by zero. So this allows to    apply Doob's optional stopping theorem. If one tries to generalize Lemma~\ref{lem:EW} to unbounded stopping times, one needs additional integrability conditions. 
Those may not be met for a small prior parameter $\alpha_S$. Therefore, we prove Lemma~\ref{lem:EW} only for bounded $\si$ and make the subsequent verification only on the interval $[S,0]$.
\end{rem}
Note that by decreasing the lower limit in the sum of line~\eqref{eq:line3} from $2$ to 1, we get for the expectation conditioned on the event $\{\Delta_t=1\}$ analogously to proof of Lemma~\ref{lem:EW}
\begin{align*}
    \IE_{t,k} & [W(\si,N_\si,\Delta_{\si})|\Delta_t=1]\notag\\&= W(t,k,1)
    +\IE_{t,k} \left[ \int_t^\si \lambda(s,N_{s-})\big(W(s,N_{s-}+1,1)-W(s,N_{s-},0) \big) \td s \middle|\Delta_t=1\right]
    \notag\\ &\quad +
    \IE_{t,k} \left[ \int_t^\si \partial_t W(s,N_s,0) \td s\middle|\Delta_t=1\right]
    \notag\\ &\quad +
    \IE_{t,k} \left[
    \sum_{1\leq j\leq k^*} \big\{ W(\sigma_t^{j-1},N_{\sigma_t^{j-1}},0)-W(\sigma_t^{j-1},N_{\sigma_t^{j-1}},1) \big\}\middle|\Delta_t=1\right].
\end{align*}

\begin{proof}[Proof of Theorem~\ref{thm:valueverifi}]

We verify that the solution of the System~\eqref{eq:HJB} is actually the value function.
To this end, let $w$ be a solution of Equation~\eqref{hjb for 0 comp} with boundary condition~\eqref{boundary condi}.
By setting $W(t,k,0):=w(t,k)$ and $W(t,k,1):=\min\{c(t,k),W(t,k,0)\}$ on $[S,0]$ we get by Theorem~\ref{lem:uniquePDE} a bounded function that satisfies  Equations~\eqref{eq:HJB}.
Let $k^*,\si, $ and $\sigma_t^j,j\in\{0,...,k^*\}$ as in Lemma~\ref{lem:sum1}.
From Lemma~\ref{lem:EW}, and since $W$ fulfills Equation~\eqref{eq:HJB} on $[S,0]$, we obtain
\begin{align}
    \IE_{t,k}  W(\si,N_{\si},\Delta_\si)
    &= W(t,k,0) +
    \IE_{t,k} 
    \sum_{2\leq j\leq k^*} \big\{ W(\sigma_t^{j-1},N_{\sigma_t^{j-1}},0)-W(\sigma_t^{j-1},N_{\sigma_t^{j-1}},1) \big\}\label{eq:termsNeed0}
    \\
    &\geq  W(t,k,0)\label{eq:V>W}
\end{align}
and analogously
\begin{align}
    \IE_{t,k} [ W(\si,N_{\si},\Delta_{\sigma})|\Delta_t=1]
    &\geq  W(t,k,1).\label{eq:V>W1}
\end{align}
Putting the result in relation to the value function, we obtain for all $\tau\in \cT_t$ with $\tau\leq\tau_0^1$:
\begin{align}\label{eq:18veri}
    \IE_{t,k}[c(\tau,N_\tau)]&=
        \IE_{t,k}[c(\tau,N_\tau)\sI_{\{\tau\leq0\}}] +
        \IE_{t,k}[c(\tau,N_\tau)\sI_{\{\tau>0\}}]
        \\
        &\stackrel{\eqref{eq:HJB}}\geq 
        \IE_{t,k}[W(\tau,N_\tau,1)\sI_{\{\tau\leq0\}}] +
        \IE_{t,k}[
            \IE_{0,N_0} [c(\tau,N_\tau)\sI_{\{\tau>0\}}]]\label{eq:verification2}\\
        &=    
        \IE_{t,k}[W(\tau,N_\tau,1)\sI_{\{\tau\leq0\}}] +
        \IE_{t,k}[
            \IE_{0,N_0} [c(\tau,N_\tau)]\sI_{\{\tau>0\}}]\\
        &\stackrel{\tiny \makebox[0pt]{{\ref{C3}}}}\geq    
        \IE_{t,k}[W(\tau,N_\tau,1)\sI_{\{\tau\leq0\}}] +
        \IE_{t,k}[
            W(0,N_0,0)\sI_{\{\tau>0\}}]\label{eq:where C1 is used}
            \\
        &=    
        \IE_{t,k}[W(\tau\wedge0,N_{\tau\wedge0},\Delta_{\tau\wedge0})].\label{eq:22veri}
\end{align}
Equations~(\ref{eq:18veri} -- \ref{eq:22veri}) also hold when conditioning on $\{\Delta_t=1\}$.
Together with Equations~\eqref{eq:V>W} and \eqref{eq:V>W1} we obtain,
$\IE_{t,k}[c(\tau,N_\tau)|\Delta_t=\Delta]\geq\IE_{t,k}[ W(\tau\wedge 0, N_{\tau\wedge 0},\Delta_{\tau\wedge0})|\Delta_t=\Delta]\geq  W(t,k,\Delta)$ for all $\tau\in \cT_t$.
So  $V(t,k,\Delta)=\inf_\tau (\IE_{t,k} [c(\tau,N_\tau)|\Delta_t=\Delta])\geq W(t,k,\Delta)$.

Let $\tau=\tau^*$ be defined as 
\eqref{optiStop} and $\si=\tau\wedge 0$. 
Then
the terms $W(\tau_t^{k-1},N_{\tau_t^{k-1}},0)$ and $W(\tau_t^{k-1},N_{\tau_t^{k-1}},1)$ from Equation~\eqref{eq:termsNeed0} are equal if $k-1<k^*$. As a result Inequalities~\eqref{eq:V>W} and~\eqref{eq:V>W1}  become equalities for $\tau^*$. 
Having $\tau=\tau^*$, stop only if $c(s,k)\leq W(s,k,0)$, Inequalities~\eqref{eq:verification2} and~\eqref{eq:where C1 is used} become equalities, too. 
Hence 
\begin{align*}
    V(t,k,\Delta)\leq\IE_{t,k} [c(\tau^*,N_{\tau^*})|\Delta_t=\Delta]&=
    W(t,k,\Delta).
\end{align*}
All in all: 
\begin{align*}
W(t,k,\Delta)\leq\inf _{\tau\in \cT_t}(\IE_{t,k} [c(\tau,N_\tau)|\Delta_t=\Delta]) =    V(t,k,\Delta)\leq\IE_{t,k} [c(\tau^*,N_{\tau^*})|\Delta_t=\Delta]=
    W(t,k,\Delta).
\end{align*}
So $W=V$ on $[S,0]$. 
\end{proof}
Note that since $V$ is unique, the solution of Equation~\eqref{hjb for 0 comp} is also unique.

\section{A Comparison with the 'One-Step-Look-Ahead' Rule }

In this section, we analyze the 'One-Step-Look-Ahead' (OSLA) rule and compare it to the optimal stopping rule from the previous sections.

The OSLA rule describes the strategy to stop at the first free parking lot at which the cost from stopping immediately is at least as small as the expected cost from continuing to the next free parking lot and stopping there. More precisely, the OSLA rule stops at the first time $t$ at which the driver enters a state $(t,k,1)$ for which  \begin{align}\label{eq:OSLAcrit}
    c(t,k)\le \IE_{t,k}[c(\tau_t^1,k+1)].
\end{align}

It is a well-known result that the OSLA rule is optimal if the stopping problem is monotone (see, e.g., \cite{ferguson2018optimal}). A stopping problem is monotone if, whenever the OSLA criterion~\eqref{eq:OSLAcrit} is satisfied at some stopping opportunity, it is almost surely satisfied at the next stopping opportunity and hence all subsequent ones.  However, our stopping problem is not monotone in general, as can be seen in Remark \ref{counter_OSLA} after Proposition \ref{OSLA_prop} below.

Note that monotonicity is only a sufficient criterion for the optimality of the OSLA rule.
Since we have a non-monotone problem, we conjecture that the OSLA rule and the globally optimal rule derived from the previous sections in general do not coincide.

In this section we demonstrate numerically that the OSLA rule cannot be generally optimal.
We do so for the specific cost function from Example~\ref{example_cost_absolute}, namely \(c(t,k)=|t|\). Here the stopping region for the OSLA rule can be determined analytically. 
We remark that in \cite[Chapter 2]{Tamaki_1985}  it is already shown that in the discrete-time version of the model the OSLA rule is in general not optimal.

\begin{prop}
\label{OSLA_prop}
    The OSLA rule for the cost function \(c(t,k)=|t|\) is characterized by: 
\begin{center}

   \fbox{
   Stop at the first time a driver enters a state $(t,k,1)$ in which \(t\ge \beta_0(2^{-\frac{1}{\alpha_S +k-1}}-1).\)}

\end{center}
    
\end{prop}
\begin{proof}
Let  a driver be in state $(t,k,0)$. 
If $t\geq 0$, Equation~\eqref{eq:OSLAcrit} is fulfilled for all subsequent free parking lots due to Condition~\ref{C3}.
Here, the stopping criterion is also fulfilled:  \(\beta_0(2^{-\frac{1}{\alpha_S +k-1}}-1) < 0 \le t\).

Let $t<0$.
We denote the expected costs of the next free parking lot after \(t\) by \(\tilde{C}(t,k):=\IE_{t,k}[c(\tau_{t}^1, k+1)]=\IE_{t,k}[|\tau_t^1|]\). Note that the Markov property of $N$ implies \(\tilde{C}(t,k)  
=\IE_{t,k}[|\tau_t^1| \mathbf{1}_{\{\tau_t^1\le 0\}}+ R(k)\mathbf{1}_{\{\tau_t^1> 0\}}]
,\) 
where \(R(k)\) is the boundary condition, i.e. \(R(k)=\IE_{0,k}  [|\tau_0^1|]=\frac{\beta_0}{\alpha_S +k-1}\).

We calculate now explicitly the value of \(\tilde{C}(t,k)\). First note that, with \eqref{surv prob}, 
\begin{align*}
  \IE_{t,k}[R(k)\mathbf{1}_{\{\tau_t^1> 0\}}]&=R(k) P_{t,k}(\tau_t^1>0)=R(k)\left(\frac{\beta_t}{\beta_0}\right)^{\alpha_S +k}, 
\end{align*}
and secondly,
\begin{align*}
    \IE_{t,k}[|\tau_t^1| \mathbf{1}_{\{\tau_t^1\le 0\}}]&=\int_t^0 |u| \frac{\alpha_S +k}{\beta_u}\left(\frac{\beta_t}{\beta_u}\right)^{\alpha_S +k} \, du\\
    &=(\alpha_S +k) \beta_t^{\alpha_S +k} \int_t^0 (-u) \left(\frac{1}{\beta_u}\right)^{\alpha_S +k+1}\, du.
\end{align*}
Substitute now \(x = \beta_u\). Note that \(u=x-\beta_S +S=-\beta_0+x\). We get
\begin{align*}
   &  (\alpha_S +k) \beta_t^{\alpha_S +k} \int_t^0 (-u) \left(\frac{1}{\beta_u}\right)^{\alpha_S +k+1}\td u \\
   &=(\alpha_S +k) \beta_t^{\alpha_S +k}\int_{\beta_t}^{\beta_0} (\beta_0-x) x^{-(\alpha_S +k+1)} \, dx\\
    &= (\alpha_S +k) \beta_t^{\alpha_S +k} \left[\beta_0\int_{\beta_t}^{\beta_0} x^{-(\alpha_S +k+1)} \, dx - \int_{\beta_t}^{\beta_0}  x^{-(\alpha_S +k)} \, dx\right] \\
    &=(\alpha_S +k) \beta_t^{\alpha_S +k} \Bigg[\frac{\beta_0}{\alpha_S +k} \left(\beta_t^{-(\alpha_S +k)}-\beta_0^{-(\alpha_S +k)}\right) -\\
    & \hspace{1cm} - \frac{1}{\alpha_S +k-1}\left(\beta_t^{-(\alpha_S +k-1)}-\beta_0^{-(\alpha_S +k-1)}\right)\Bigg] \\
    &= \beta_0 - \beta_0 \left(\frac{\beta_t}{\beta_0}\right)^{\alpha_S +k} - \frac{\alpha_S +k}{\alpha_S +k-1} \left(\beta_t-\beta_0 \left(\frac{\beta_t}{\beta_0}\right)^{\alpha_S +k} \right).
\end{align*}

Therefore
\begin{align*}
    \tilde{C}(t,k)&= \beta_0 - \beta_0 \left(\frac{\beta_t}{\beta_0}\right)^{\alpha_S +k} - \frac{\alpha_S +k}{\alpha_S +k-1} \left(\beta_t-\beta_0 \left(\frac{\beta_t}{\beta_0}\right)^{\alpha_S +k} \right) + R(k)\left(\frac{\beta_t}{\beta_0}\right)^{\alpha_S +k}\\
    &=\beta_0 - \beta_0 \left(1+\frac{t}{\beta_0}\right)^{\alpha_S +k} - \frac{\alpha_S +k}{\alpha_S +k-1} \left(\beta_t-\beta_0 \left(1+\frac{t}{\beta_0}\right)^{\alpha_S +k} \right)\\
    & \quad + R(k)\left(1+\frac{t}{\beta_0}\right)^{\alpha_S +k}.
\end{align*}

Now we characterize the stopping region for a fixed \(k\). In particular, we  find an explicit expression for \(S_k:=\{t \in [S,0]: |t| \le \tilde C (t,k)\}.\) 
Calculations show
\begin{align*}
-t&\le \tilde C(t,k)= \beta_0 - \beta_0 \left(1+\frac{t}{\beta_0}\right)^{\alpha_S +k} - \frac{\alpha_S +k}{\alpha_S +k-1} \left(\beta_t-\beta_0 \left(1+\frac{t}{\beta_0}\right)^{\alpha_S +k} \right) \\
&\quad+ R(k)\left(1+\frac{t}{\beta_0}\right)^{\alpha_S +k}.   
\end{align*}
We substitute \(x=1+\frac{t}{\beta_0}.\) Note that \(t=\beta_0(x-1).\) It follows
\begin{align*}
    \beta_0(1-x)&\le \beta_0 (1-x^{\alpha_S +k}) - \frac{\alpha_S +k}{\alpha_S +k-1} \left[(t+\beta_0)-\beta_0x^{\alpha_S +k}\right] + R(k) x^{\alpha_S +k}.
\end{align*}
This is equivalent to
\begin{align*}
    \beta_0(1-x) - \beta_0 (1-x^{\alpha_S +k})&\le - \frac{\alpha_S +k}{\alpha_S +k-1}\beta_0 \left[x-x^{\alpha_S +k}\right] + R(k) x^{\alpha_S +k},
\end{align*}
and to
\begin{align*}
  &(x^{\alpha_S +k}-x)\beta_0 \left(1-\frac{\alpha_S +k}{\alpha_S +k-1}\right) \le R(k) x^{\alpha_S +k}. 
\end{align*}
Solving for \(x\) yields
\begin{align*}
    x\ge\left(\frac{\beta_0}{\beta_0+R(k)(\alpha_S +k-1)}\right)^{\frac{1}{\alpha_S +k-1}}.
\end{align*}
Lastly,
\begin{align*}
    t\ge \beta_0\left[\left(\frac{\beta_0}{\beta_0+R(k)(\alpha_S +k-1)}\right)^{\frac{1}{\alpha_S +k-1}}-1\right].
\end{align*}

With \(R(k)=\frac{\beta_0}{\alpha_S +k-1}\), the expression simplifies to 
\begin{align*}
    t\ge \beta_0\left[\left(\frac12\right)^{\frac{1}{\alpha_S +k-1}}-1\right]=\beta_0(2^{-\frac{1}{\alpha_S +k-1}}-1).
\end{align*}

We conclude that \(S_k=[\max\{S,t_k^{OSLA}\},0],\) where \(t_k^{OSLA}=\beta_0(2^{-\frac{1}{\alpha_S +k-1}}-1).\)

\end{proof}

\begin{rem}
\label{counter_OSLA}
Suppose that \(S=-1, \alpha_S=\beta_S=2\). The event \(\{\tau_S^1 \in (-0.8, -0.7), \tau_S^2\in (-0.69,-0.64)\}\)  has a strictly positive probability.  Then for \(c(t,k)=|t|\), according to Proposition \ref{OSLA_prop}, the OSLA criterion is satisfied at \(\tau_S^1\) if it is greater than \(-0.87\) and it is not satisfied at \(\tau_S^2\) if \(\tau_S^2<-0.62\). Hence, the OSLA criterion on the event \(\{\tau_S^1 \in (-0.8, -0.7), \tau_S^2\in (-0.69,-0.64)\}\) is satisfied for \(\tau_S^1\), but not for \(\tau_S^2\).  Therefore, our stopping problem is not monotone. 
\end{rem}
We now present results from numerical experiments that illustrate the difference between the globally optimal stopping rule and the OSLA rule. Table \ref{tab:results} depicts stopping regions of both rules for four different sets of initial parameters $\alpha_S$ and $\beta_S$. We use the following notation: \(t^{OSLA}_k\) is the first time the driver would stop if there were a free parking lot, provided she is using the OSLA rule and has observed \(k\) free parking lots including the current free parking spot. Furthermore, \(t^*_k\) is defined in exactly the analogous way, except that the globally optimal stopping rule is used instead of the OSLA rule.

Since we do not find a closed-form solution of the system \eqref{eq:HJB}, we solve it numerically. 
We discretize \([S,0]\) by the uniform grid \(t_i=S+ih, \, i= 0, \dots, M, \, h=|S|/M,\) and approximate the value function by a backward explicit Euler scheme, starting from the terminal condition \(W(0,k,0)=\frac{\beta_S -S}{\alpha_S+k-1}\) and using
\begin{align*}
    W(t_i,k,0)&=W(t_{i+1}, k ,0) - h \lambda(t_{i+1}, k) \left(W(t_{i+1},k,0)-W(t_{i+1},k+1,1)\right)
\intertext{and}
     W(t_i,k,1)&=\min\{|t_i|,\,W(t_i,k,0)\}.
\end{align*}

Since the system is coupled upwards in \(k\), we truncate it at a level \(\tilde K\) and impose the artificial boundary condition \(W(t,\tilde K,1)=0\) for \(t \in [S,0].\)

The mesh size \(h\) was chosen as \(h=\frac{1}{1000}\) and the truncation level as
\(\tilde K=1000\). Further
refinement of the mesh or an increase of \(\tilde K\) does not lead to a noticeable change
in the computed stopping regions $\{(t,k): |t|\le W(t,k,0)\}$.

\begin{table}[H]
  \centering
  \caption{Numerical results for different parameter pairs $(\alpha_S ,\beta_S)$ and $c(t,k)=|t|$.} 
  \label{tab:results}
  \begin{tabular}{r ccc ccc}
    \toprule
    & \multicolumn{3}{c}{$S=-1,\,  \alpha_S =2,\ \beta_S=2$}
    & \multicolumn{3}{c}{$S=-1,\  \alpha_S =10,\ \beta_S=10$} \\
    \cmidrule(lr){2-4} \cmidrule(lr){5-7}
    $k$ & $t^{OSLA}_k$ & $t^*_k$ & $t^*_k-t^{OSLA}_k$
        & $t^{OSLA}_k$ & $t^*_k$ & $t^*_k- t^{OSLA}_k$ \\
    \midrule
    1  & $-0.879$ & $-0.837$ & $0.042$ & $-0.737$ & $-0.734$ & $0.003$ \\
    2  & $-0.619$ & $-0.602$ & $0.017$ & $-0.672$ & $-0.669$ & $0.003$ \\
    3  & $-0.477$ & $-0.469$ & $0.008$ & $-0.617$ & $-0.615$ & $0.002$ \\
    4  & $-0.388$ & $-0.383$ & $0.005$ & $-0.571$ & $-0.569$ & $0.002$ \\
    5  & $-0.327$ & $-0.324$ & $0.003$ & $-0.531$ & $-0.530$ & $0.001$ \\
    \addlinespace
    10 & $-0.183$ & $-0.182$ & $0.001$ & $-0.394$ & $-0.393$ & $0.001$ \\
    20 & $-0.097$ & $-0.096$ & $0.001$ & $-0.260$ & $-0.259$ & $0.001$ \\
    \midrule
    & \multicolumn{3}{c}{$S=-1,\ \alpha_S =10,\ \beta_S=5$}
    & \multicolumn{3}{c}{$S=-1,\ \alpha_S =5,\ \beta_S=10$} \\
    \cmidrule(lr){2-4} \cmidrule(lr){5-7}
    $k$ & $t^{OSLA}_k$ & $t^*_k$ & $t^*_k-t^{OSLA}_k$
        & $t^{OSLA}_k$ & $t^*_k$ & $t^*_k - t^{OSLA}_k$ \\
    \midrule
    1  & $-0.402$ & $-0.400$ & $0.002$ & $-1.424 \vee -1$ & $-1$ & $0$ \\
    2  & $-0.366$ & $-0.365$ & $0.001$ & $-1.200\vee -1$ & $-1$ & $0$ \\
    3  & $-0.337$ & $-0.335$ & $0.002$ & $-1.037\vee -1$ & $-1$ & $0$ \\
    4  & $-0.312$ & $-0.310$ & $0.002$ & $-0.913$ & $-0.908$ & $0.005$ \\
    5  & $-0.290$ & $-0.289$ & $0.001$ & $-0.815$ & $-0.811$ & $0.004$ \\
    \addlinespace
    10 & $-0.215$ & $-0.214$ & $0.001$ & $-0.531$ & $-0.530$ & $0.001$ \\
    20 & $-0.142$ & $-0.141$ & $0.001$ & $-0.313$ & $-0.312$ & $0.001$ \\
    \bottomrule
  \end{tabular}
\end{table}

The calculations from Table \ref{tab:results} are in line with our intuition derived from the missing monotonicity. The stopping region of the OSLA rule is strictly larger than that of the optimal stopping rule except for \(\alpha_S=5, \beta_S=10\) and \(k\le 3\), where both regions cover the whole street. 

The expected costs of the next free parking lot are greater than or equal to the continuation value.  
As a result, under the OSLA rule there is a tendency to stop too early.  
Furthermore, Table~\ref{tab:results} indicates that for large \(\alpha_S+k\) the OSLA rule seems almost optimal. The differences are largest when \(\alpha_S+k\) is small and \(\beta_S\) is not too large compared to the length \(|S|\) of the street. The jump times are essential for calculating optimal strategies and they cause a break in the monotonicity of the problem.  A small value of \(\alpha_S+k\), with small \(\beta_S\), corresponds to a driver with much uncertainty about the jump intensity, and each newly observed free parking lot changes her estimate considerably. Especially early jumps, i.e. early free parking lots, may have a very strong influence on subsequent actions, since the posterior probability changes substantially. 
However, the OSLA strategy ignores that the estimated jump intensity will increase at future jump times and therefore may stop too early if there is a possibility of parking. If, on the other hand, \(\beta_S\) is large compared to \(|S|\), only a few free parking lots are expected on the street, and both rules nearly coincide. 

\bibliographystyle{abbrvnat}
\bibliography{literatur,lit}

\end{document}